# Random perturbations of codimension one homoclinic tangencies in dimension 3

Vítor Araújo (`vdaraujo@fc.up.pt`) [*]
Centro de Matemática da Universidade do Porto
4099-002 Porto, Portugal
Telefone: +351 22340 1461; Fax: +351 22340 1453

**Abstract**

Adding small random parametric noise to an arc of diffeomophisms of a manifold of dimension 3, generically unfolding a codimension one quadratic homoclinic tangency $q$ associated to a sectionally dissipative saddle fixed point $p$, we obtain not more than a finite number of *physical* probability measures, whose ergodic basins cover the orbits which are recurrent to a neighborhood of the tangency point $q$. This result is in contrast to the extension of Newhouse's phenomenon of coexistence of infinitely many sinks obtained by Palis and Viana in this setting.

There is a similar result for the simpler bidimensional case whose proof relies on geometric arguments. We now extend the arguments to cover three dimensional manifolds.

## 1 Introduction

Let $f: M \to M$ be a smooth ($C^\infty$) diffeomorphism of a compact manifold with dimension 3 having a hyperbolic fixed point $p$ of saddle type with index 2. Let us assume that the stable and unstable manifolds of $p$ have an homoclinic tangency at a point $q$ and that $p$ is *sectionally dissipative*: if $|\sigma| > 1 > |\lambda_1| \geq |\lambda_2|$ are the eigenvalues of $df_p$, then $|\sigma\lambda_i| < 1, i = 1, 2$.

In [10] Palis and Viana proved that near such map $f$ (in the $C^2$ topology) there is a residual subset of an open set of diffeomorphisms such that each of its elements exhibits infinitely many coexisting sinks. In particular a generic unfolding of such a (quadratic) homoclinic tangency by a one-parameter family of diffeomorphisms admits residual subsets of intervals in the parameter line whose corresponding diffeomorphisms display infinitely many sinks. Moreover those intervals accumulate on the parameter of tangency.

---

[*]Work supported by *Fundação para a Ciência e a Tecnologia* (FCT) through *Centro de Matemática da Universidade do Porto* (CMUP) and grant SAPIENS/FCT 36581/99. Available as a PDF file from `http://www.fc.up.pt/cmup`.

Here we prove that, by randomly perturbing a parameter close to the parameter of tangency — not excluding parameters corresponding to infinitely many sinks — and considering only those orbits which are recurrent to a given neighborhood of a fixed tangency point, their mean sojourn times (Birkhoff averages) are given by a finite number of absolutely continuous stationary probability measures.

This is an extension to a higher dimension of the geometric proof of Theorem 2 in [2] for a one-parameter family of diffeomorphisms generically unfolding a quadratic homoclinic tangency associated to a dissipative saddle fixed point on a surface.

Very little is known about this generalization of *Newhouse's phenomenon* to higher dimensions, even about its original version for surface diffeomorphisms (see [9] and also [5] for more details), apart from the fact that it exists. A very interesting question is whether in the parameter line this is a rare phenomenon with respect to Lebesgue measure — still unanswered and apparently rather difficult to attack. Another question concerning the relevance of this phenomenon is the question of stochastic stability: is the asymptotic behavior of maps displaying infinitely many coexisting sinks stable under random perturbations?

The main result is a contribution to an answer to the second question, since we construct the analogues of $SRB$ measures in the non-random setting. See also [3] for sufficient conditions ensuring the convergence of these $SRB$ measures, i.e., the stochastic stability of these maps.

## 1.1 Random perturbations

We are interested in studying random perturbations of the map $f$ along an arc of diffeomorphisms. For that, we take a one-parameter family $\Phi : ]-1,1[ \to C^2(M,M)$, $t \mapsto f_t$, with $f = f_0$. We assume that this family generically unfolds the quadratic homoclinic tangency when $t$ grows past 0, see the statement of the main theorem.

Given $x \in M$ we call the sequence $\left(f_{\underline{t}}^n x\right)_{n \geq 1}$ a *random orbit* of $x$, where $\underline{t}$ denotes an element $(t_1, t_2, t_3, \dots)$ in the product space $T = I^{\mathbb{N}}$, where $I = ]-1,1[$ and $f_{\underline{t}}^n = f_{t_n} \circ \cdots \circ f_{t_1}$ for $n \geq 1$.

We also take a family $(\theta_\epsilon)_{\epsilon > 0}$ of probability measures on $I$ such that $(\operatorname{supp} \theta_\epsilon)_{\epsilon > 0}$ is a nested family of connected compact sets and $\operatorname{supp} \theta_\epsilon \to \{0\}$ when $\epsilon \to 0$. We will refer to $\{\Phi, (\theta_\epsilon)_{\epsilon > 0}\}$ as a *random perturbation* of $f$.

## 1.2 Physical probability measures

In what follows *Lebesgue measure* will be denoted by $m$ and will mean a normalized Riemannian volume form on the compact manifold $M$ fixed for the rest of the paper. In the setting of random perturbations of a map we say that a Borel probability measure $\mu^\epsilon$ on $M$ is *physical* if for a positive Lebesgue measure set of points $x \in M$, the averaged sequence of Dirac probability measures $\delta_{f_{\underline{t}}^n(x)}$ along random orbits $\left(f_{\underline{t}}^n x\right)_{n \geq 0}$ converges in the weak$^*$ topology to $\mu^\epsilon$ for $\theta_\epsilon^\mathbb{N}$ almost every $\underline{t} \in T$. That is,

$$\lim_{n \to +\infty} \frac{1}{n} \sum_{j=0}^{n-1} \varphi(f_{\underline{t}}^n x) = \int \varphi \, d\mu^\epsilon \quad \text{for all continuous } \varphi \colon M \to \mathbb{R} \tag{1}$$

and $\theta_\epsilon^\mathbb{N}$ almost every $\underline{t} \in T$. We denote the set of points $x \in M$ for which (1) holds by $B(\mu^\epsilon)$ and call it the *basin of $\mu^\epsilon$*.

## 1.3 Main result

The family of $C^\infty$ diffeomorphisms $\Phi \colon ]-1,1[ \to C^\infty(M,M)$ of a compact manifold $M$ with dimension 3 that we will be considering satisfies the following conditions, see figure 1.

1. The map $f = \Phi(0)$ admits a hyperbolic saddle fixed point $p$ such that:

    (a) $df_p$ has a single expanding eigenvalue $\sigma > 1$;
    (b) for every other eigenvalue $\lambda$ it holds $|\sigma \lambda| < 1$;
    (c) $df_p$ admits a *least contracting eigenvalue* $\lambda_1$, satisfying $|\lambda_1 > |\lambda_2|$;
    (d) the invariant manifolds $W_f^u(p)$ and $W_f^s(p)$ have a quadratic tangency at a point $q$;
    (e) there is a compact neighborhood $U$ whose maximal invariant subset $\Lambda_p = \cap_{n \in \mathbb{Z}} f^n(U)$ is a hyperbolic basic set containing $p$.

2. In a neighborhood of $(0,p)$ in $]-1,1[ \times M$ the map $F(t,x) = \Phi(0)x$ satisfies the conditions for the existence of $C^2$ linearizing coordinates, see [13, 14], and in these coordinates

    (a) the stable manifold $W_f^s(p)$ is identified with the "box" $\{0\} \times [-2,2]^2$ and the unstable manifold $W_f^u(p)$ with the line $[-2,2] \times \{0\}$;
    (b) the tangency point $q$ has coordinates $(0,1,1)$ and there is $N \geq 1$ such that $r = f^{-N}q$ has coordinates $(1,0,0)$;
    (c) the tangent vector $B$ to $W_f^u(p)$ at $q$ ($B = D_x(f_t^N(x))\mid_{x=r}$) is transverse both to the strong-contracting and weak-contracting directions of $df_p$;
    (d) the *unfolding direction* $A = D_t(\Phi(t)(q))\mid_{t=0}$ and the weak-contracting direction $D$ together with $B$ are linearly independent, i.e. $\det(A,B,D) \neq 0$;

**Remark 1.1.** *Condition (2) above is generic in the space of $C^\infty$ families of diffeomorphisms satisfying (1) and the families satisfying (1) form an open set of families.*

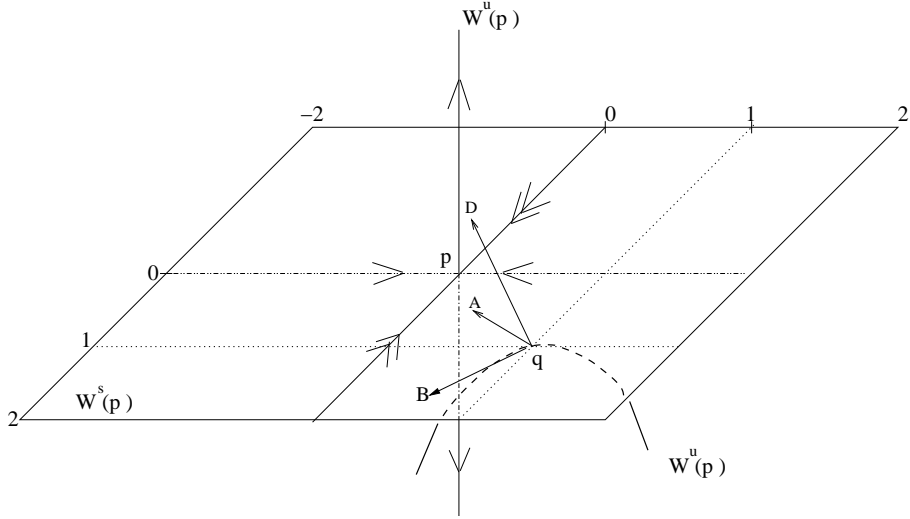

Figure 1: The directions $A, B$ and $D$ on the quadratic tangency at $q$

**Remark 1.2.** *Condition (1c) is natural in this setting since Palis and Viana [10, Section 5] have shown that, for generic one-parameter families of diffeomorphisms unfolding a codimension 1 quadratic homoclinic tangency at $t = 0$ associated to a sectionally dissipative hyperbolic fixed point $p$, there are parameter values $t_j \to 0$ and homoclinic tangencies associated to periodic points $p_j \to p$ such that $df_{t_j}^{l_j}(p_j)$ has a unique weakest contracting eigenvalue (where $l_j$ is the period of $p_j$).*

**Remark 1.3.** *Condition 1(e) is rather natural since a hyperbolic fixed point $p$ might have more homoclinic orbits besides $(f^k q)_{k \in \mathbb{Z}}$, in which case the saddle is part of a non-trivial hyperbolic basic set $\Lambda_p \ni p$. This is the case for a denumerable set of parameters after unfolding a homoclinic tangency [9, Chpt. 3]. Nevertheless the set $\Lambda_p$ might trivially be $\{p\}$, in which case the results that follow still apply.*

To specify what orbits our result applies to we need the following notions. We will focus on the behavior of orbits that eventually stay in a *closed neighborhood $U$* of the basic set together with a neighborhood $\mathcal{Q}$ of $q$ and that pass infinitely many times through $\mathcal{Q}$.

**Definition 1.4.** *Given a point $x \in M$ and a sequence $\underline{t} \in I^{\mathbb{N}}$ we define the first return time to $\mathcal{Q}$ by $r(1, \underline{t}, x) = \min\{n \geq 1 : f_{\underline{t}}^n x \in \mathcal{Q}\}$, with the convention that $\min \emptyset = \infty$.*

*Recursively we may now define the $(k+1)$-th return time as a function of the $k$th return time $r(k+1, \underline{t}, x) = r(1, \sigma^{r(k,\underline{t},x)} \underline{t}, f_{\underline{t}}^{r(k,\underline{t},x)} x)$, where $\sigma : I^{\mathbb{N}} \to I^{\mathbb{N}}$ is the left shift on sequences.*

As explained in the Introduction, we are going to study the mean sojourn times of perturbed orbits that are recurrent to a neighborhood $\mathcal{Q}$ of the tangency point $q$. These orbits return infinitely often to $\mathcal{Q}$ under a set of perturbations of positive probability, in the following sense.

**Definition 1.5.** *Given a probability measure $\theta$ on $I$ we say that a point $x$ is $\theta$ recurrent if there is a set $W \subset I^{\mathbb{N}}$ such that*

- *$\theta^{\mathbb{N}}(W) > 0$;*
- *for every $\underline{t} \in W$ there is $n_0 = n_0(\underline{t}) \in \mathbb{N}$ such that $f_{\underline{t}}^k x \in U \cup \mathcal{Q}$ for all $k \geq n_0$;*
- *$r(k, \underline{t}, x) < \infty$ for all $k \geq 1$ and every $\underline{t} \in W$.*

We are now ready for the statement of the main result. For definiteness we set

$$\theta_\epsilon = (2\epsilon)^{-1} \text{Leb } ]]t_0 - \epsilon, t_0 + \epsilon[ \tag{2}$$

for given $t_0$ — the normalized restriction of Lebesgue measure to the $\epsilon$-neighborhood of $t_0$ on the real line.

**Theorem A.** *For every $C^\infty$ one-parameter family of diffeomorphisms of a compact manifold with dimension 3 satisfying conditions (1) and (2), and for each quadratic homoclinic point $q$ associated to the sectionally dissipative saddle $p$, there are a closed neighborhood $\mathcal{Q}$ of $q$ and $t_\star > 0$ such that, for each $t_0 \in ]0, t_\star[$ and $\epsilon \in ]0, t_\star - t_0[$, the random perturbation $\{\Phi, (\theta_\epsilon)_\epsilon\}$ admits a finite number of physical probability measures $\mu_1, \ldots, \mu_l$ whose support intersects $\mathcal{Q}$ and*

1. *$\mu_i \perp \mu_j$ for all $1 \leq i < j \leq l$ (the physical measures are mutually singular);*
2. *$\mu_i \ll m$ for all $i = 1, \ldots, l$;*
3. *for all $\theta_\epsilon$ recurrent $x \in \mathcal{Q}$ with respect to a measurable set $W$ there exist measurable sets $W_1 = W_1(x), \ldots, W_l = W_l(x) \subset W$ satisfying*

   (a) *$W_i \cap W_j = \emptyset$ for all $1 \leq i < j \leq l$;*
   (b) *$\theta_\epsilon^{\mathbb{N}}(W \setminus (W_1 \cup \cdots \cup W_l)) = 0$;*
   (c) *for every $i = 1, \ldots, l$ and all $\underline{t} \in W_i$ we have*

   $$\lim_{n \to \infty} \frac{1}{n} \sum_{j=0}^{n-1} \varphi(f_{\underline{t}}^j x) = \int \varphi \, d\mu_i \quad \text{for all} \quad \varphi \in C^0(M, \mathbb{R}).$$

The theorem assures the existence of a finite number of physical probability measures for the random perturbation $\{\Phi, (\theta_\epsilon)_\epsilon\}$ which describe the asymptotics of the Birkhoff averages of almost all perturbed orbits of every recurrent point to the neighborhood $\mathcal{Q}$. In this sense these probability measures capture the asymptotic behavior of the system near the homoclinic point.

**Remark 1.6.** *This Theorem holds true for probability measures $\theta_\epsilon$ equivalent to Lebesgue measure, i.e. $\theta_\epsilon = \phi_\epsilon \cdot m$ with $\phi_\epsilon > 0$ Lebesgue a.e. on their support, which equals $[t_0 - d(\epsilon), t_0 + d(\epsilon)]$, where $d : \mathbb{R}^+ \to \mathbb{R}^+$ is a monotonous increasing map such that $\lim_{\epsilon \to 0} d(\epsilon) = 0$.*

The proof is organized as follows: in section 2 we present some measure theoretic results used in the sequel and in Sections 3 and 4 we prove Theorem A.

In Section 3 a linearization of the coordinates near the hyperbolic saddle is used to obtain properties of the return iterates to $\mathcal{Q}$ of randomly perturbed orbits. In Section 4 properties of accumulation points of randomly perturbed orbits are obtained and then used, together with the properties from Section 3, to prove Theorem A.

The strategy is to show that fixing $\mathcal{Q}, t_0$ and $\epsilon > 0$ as in the statement of the result, then given any stationary ergodic probability measure $\mu$ for the random perturbation $\{\Phi, \theta_\epsilon\}$ whose support intersects $\mathcal{Q}$, the volume of the basin is bounded from below by a constant depending only on the family $\Phi$ and on the value of $\epsilon$. Since the basins of different probability measures are disjoint, this shows that there can only be a finite number of such measures.

In section 5 we try to explain some of the hurdles we face when attempting to extend the same strategy of proof to a higher dimensional setting (diffeomorphisms on manifolds with dimension $\geq 4$).

## 2 Stationary probability measures

In the setting of random perturbations of a map, we say that a set $A \subset M$ is *invariant* if $f_t(A) \subset A$, at least for $t \in \text{supp}(\theta_\epsilon)$ with $\epsilon > 0$ small. The usual invariance of a measure with respect to a transformation is replaced by the following one: a probability measure $\mu$ is said to be *stationary*, if for every continuous $\varphi : M \to \mathbb{R}$ it holds

$$\int \varphi \, d\mu = \iint \varphi(f_t x) \, d\mu(x) \, d\theta_\epsilon(t). \tag{3}$$

It is not difficult to see (cf. [2]) that a stationary measure $\mu$ satisfies

$$x \in \text{supp}(\mu) \quad \Rightarrow \quad f_t x \in \text{supp}(\mu) \quad \text{for all} \quad t \in \text{supp}(\theta_\epsilon) \tag{4}$$

just by continuity of $\Phi$. This means that if $\mu$ is a stationary measure, then $\text{supp}(\mu)$ is an invariant set. Let us now fix $x \in M$ and consider

$$\mu_n(x) = \frac{1}{n} \sum_{j=0}^{n-1} (f_x^j)_* \theta_\epsilon^\mathbb{N}. \tag{5}$$

Here $(f_x^j)_* \theta_\epsilon^\mathbb{N}$ is the push-forward of $\theta_\epsilon^\mathbb{N}$ to $M$ via $f_x^j : T \to M$, defined as $f_x^j(\underline{t}) = f_{\underline{t}}^j x$. Since this is a sequence of probability measures on the compact manifold $M$, then it has weak* accumulation points.

**Lemma 2.1.** *Every weak* accumulation point of* $\big(\mu_n(x)\big)_n$ *is a stationary probability measure.*

*Proof.* See [2, Sec. 7.1] or [1, Lemma 3.5]. □

This provides us with the means for easily obtaining stationary measures under random perturbations.

## 2.1 Ergodic decomposition

We now state a result that guarantees that every stationary measure can be decomposed into a combination of ergodic measures. To present the result with a simple formulation we consider the following skew-product map obtained from the random perturbation $\{\Phi, (\theta_\epsilon)_\epsilon\}$:
$$F : T \times M \to T \times M, \quad (\underline{t}, x) \mapsto (\sigma \underline{t}, f_{t_1} x).$$

**Definition 2.2.** *A stationary probability measure $\mu$ is* ergodic *if the measure $\theta_\epsilon^\mathbb{N} \times \mu$ is $F$ ergodic or, equivalently, if every invariant set $A$ has measure $\mu(A) = 0$ or $1$.*

**Proposition 2.3.** *Let $\mu$ be a stationary probability measure for the random perturbation $\{\Phi, \theta_\epsilon\}$. Then there is a family of ergodic probability measures $\mu_{\underline{t},x}$ defined for $\theta_\epsilon^\mathbb{N} \times \mu$ almost every $(\underline{t}, x)$ such that for all continuous $\varphi : M \to \mathbb{R}$ it holds*
$$\int \varphi \, d\mu = \int \left[ \int \varphi(y) \, d\mu_{\underline{t},x}(y) \right] d(\theta_\epsilon^\mathbb{N} \times \mu)(\underline{t}, x).$$

There is a very general statement and proof in [7, Appendix A.1] and [4, Sec. 1.4] but we present a short proof based on a non-random version of this result from [8].

*Proof.* Let $\varphi : M \to \mathbb{R}$ be continuous and $\mu$ a stationary probability measure as in the statement. Then it is straightforward to check that $\theta_\epsilon^\mathbb{N} \times \mu$ is $F$ invariant. By the result on ergodic decomposition for invariant measures (see [8, Chap. II, § 6]) we know that there is a family $\eta_{\underline{t},x}$ of $F$ ergodic probability measures defined for $\theta_\epsilon^\mathbb{N} \times \mu$ almost every point $(\underline{t}, x)$ such that
$$\int \varphi \circ \pi \, d(\theta_\epsilon^\mathbb{N} \times \mu) = \int \left[ \int (\varphi \circ \pi)(\underline{s}, y) \, d\eta_{\underline{t},x}(\underline{s}, y) \right] d(\theta_\epsilon^\mathbb{N} \times \mu)(\underline{t}, x), \qquad (6)$$

where $\pi : T \times M \to M$ is the projection on the second coordinate, so that $\varphi \circ \pi$ is a continuous function on $T \times M$.

We have that $\eta_{\underline{t},x} = \theta_\epsilon^\mathbb{N} \times \mu_{\underline{t},x}$ for some probability measure $\mu_{\underline{t},x}$ which will be $F$ ergodic by construction. Indeed, by construction of $\eta_{\underline{t},x}$ we have the weak$^*$ limit
$$\tilde{\pi}_\star(\eta_{\underline{t},x}) = \lim_{n \to +\infty} \frac{1}{n} \sum_{j=0}^{n-1} \tilde{\pi}_\star(\delta_{F^j(\underline{t},x)}) = \lim_{n \to +\infty} \frac{1}{n} \sum_{j=0}^{n-1} \delta_{\sigma^j \underline{t}} = \theta_\epsilon^\mathbb{N} \quad \text{for} \quad \theta_\epsilon^\mathbb{N} \text{ a.e. } \underline{t},$$

not depending on $(\underline{t}, x)$ because $\theta_\epsilon^\mathbb{N}$ is $\sigma$ ergodic, where $\tilde{\pi} : T \times M \to T$ is the projection on the first coordinate. Finally since $\int \varphi \circ \pi \, d(\theta_\epsilon^\mathbb{N} \times \mu) = \int \varphi \, d\mu$ and $\int \varphi \circ \pi \, d\eta_{\underline{t},x} = \int \varphi \, d\mu_{\underline{t},x}$, we see that (6) equals $\int \varphi \, d\mu = \int \left[ \int \varphi \, d\mu_{\underline{t},x} \right] d(\theta_\epsilon^\mathbb{N} \times \mu)(\underline{t}, x)$ as stated. □

## 2.2 Generic points of a shift invariant product measure

The next two measure theoretic results will be essential in the arguments that follow — their proofs may be found in [2, Section 5, Lemmas 5.2 and 5.3]. Let $\theta$ be a Borel probability measure on $I$ and $d(\underline{t}, \underline{s}) = \sum_{i \geq 1} 2^{-i} |t_i - s_i|$ be the usual distance in $T$.

**Lemma 2.4.** *Let $V, W \subset T$ be such that $\theta^{\mathbb{N}}(V) \cdot \theta^{\mathbb{N}}(W) > 0$. Then for $\theta^{\mathbb{N}}$ almost every $\underline{t} \in W$ there is $k_0 \in \mathbb{N}$ such that for all $k \geq k_0$ and every $\eta > 0$*

$$\theta^{\mathbb{N}}\left(\{\underline{s} \in W : d(\underline{s}, \underline{t}) < \eta \text{ and } \sigma^k \underline{s} \in V\}\right) > 0.$$

This lemma essentially says that the elements of a positive measure set $W$ for an infinite product probability $\theta^{\mathbb{N}}$ have coordinates very well spread over the base measure space, as long as we look at big enough indexes $k$ and ignore subsets of null measure.

A slight generalization of Lemma 2.4 will be needed, dealing with *double sections* and arriving at a similar conclusion with respect to a single coordinate.

**Definition 2.5.** *Given $W \subset T$ and $\underline{t}, \underline{s} \in T$ we define a* double section *of $W$ through $\underline{t}$ and $\underline{s}$ at $k \geq 1$ by $W(\underline{t}, k, \underline{s}) = \{\underline{u} \in W : u_1 = t_1, \ldots, u_k = t_k \text{ and } u_{k+2} = s_1, u_{k+3} = s_2, \ldots\}.$*

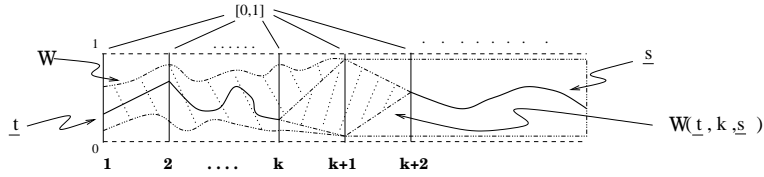

Figure 2: *Representation of the double section $W(\underline{t}, k, \underline{s})$*

If $\theta^{\mathbb{N}}(W) = 1$, then Fubini's Theorem guarantees that

$$1 = \theta^{\mathbb{N}}(W) = \int \theta^{\mathbb{N}}(\sigma^k W(\underline{t}, k)) \, d\theta^k(\underline{t}) = \iint \theta(p_{k+1} W(\underline{t}, k, \underline{s})) \, d\theta^{\mathbb{N}}(\underline{s}) \, d\theta^k(\underline{t})$$

and since $\theta$ is a probability measure, we conclude that $\theta(p_{k+1} W(\underline{t}, k, \underline{s})) = 1$ for $\theta^k$ almost all $\underline{t} \in I^k$ and $\theta^{\mathbb{N}}$ almost every $\underline{s} \in T$, where $p_{k+1}$ is the projection on the $(k+1)$th coordinate. We may think of next lemma as an extension of this result to the setting of positive measure subsets.

**Remark 2.6.** *Assume again that $\theta^{\mathbb{N}}(W) = 1$ and define for $k < l$ and $\underline{t}, \underline{s}, \underline{u} \in T$ the triple section $W(\underline{t}, k, \underline{s}, l, \underline{u})$ to be the set*

$$\{\underline{v} \in T : (v_1, \ldots, v_k) = (t_1, \ldots, t_k), (v_{k+2}, \ldots, v_l) = (s_1, \ldots, s_{l-k-1}), \sigma^{l+2}\underline{v} = \underline{u}\}.$$

*Then Fubini's Theorem will likewise ensure that $\theta^2(p_{k+1,l+1}(W(\underline{t}, k, \underline{s}, l, \underline{u}))) = 1$ for $\theta^{\mathbb{N}}$ almost all $\underline{t}, \underline{s}, \underline{u}$ and every pair of integers $k < l$, where $p_{k+1,l+1}$ is the projection on the $(k+1)$th and $(l+1)$th coordinates.*

**Lemma 2.7.** *Let $W \subset T$ be such that $\theta^{\mathbb{N}}(W) > 0$. Then for $\theta^{\mathbb{N}}$ almost every $\underline{t} \in W$ and for every $\gamma, \delta \in (0,1)$ there exists $k_0 \in \mathbb{N}$ such that for all $k \geq k_0$ there are sets $W_k \subset W$ and $Y_k \subset T$ with the properties*

1. $\underline{t} \in W_k$;

2. $\theta^{\mathbb{N}}(W_k) > 0$;

3. $\theta^{\mathbb{N}}(Y_k) \geq 1 - \gamma$;

4. $\theta\left(p_{k+1}W_k(\underline{t}, k, \underline{s})\right) \geq 1 - \delta$ for $\theta^{\mathbb{N}}$ almost all $\underline{t} \in W_k$ and $\underline{s} \in Y_k$.

That is, for generic vectors of a positive $\theta^{\mathbb{N}}$ measure subset, the measure of the double section may be made very close to that of the entire space, as long as the coordinate index is sufficiently big.

## 3 The geometric setting

In this section we prove some geometric results to be used in the proof of the Main Theorem.

### 3.1 Adapting the linearized coordinates

Let $\varphi_t : L \to \mathbb{R} \times \mathbb{R}^2$ be the mapping from a neighborhood $L$ of $p$ in $M$ that linearizes the action of $f_t = \Phi(t)$ for small $|t|$,

$$f_t \varphi_t(z, X) = \varphi_t(\sigma_t z, \Lambda_t X), \tag{7}$$

where $\sigma_t \in \mathbb{R}, |\sigma_t| > 1$ and $\Lambda_t \in \mathcal{L}(\mathbb{R}^2, \mathbb{R}^2), \|\Lambda_t\| < 1$ for some norm $\|\cdot\|$ in $\mathbb{R}^2$. By appropriately rescaling the $z$-axis and choosing the basis in $\{0\} \times \mathbb{R}^2$, we may ensure that condition (2b) holds and

- the analytic continuation $p_t$ of the saddle fixed point $p$ for $f = f_0$ satisfies $p_t = (0,0,0)$ in the $\varphi_t$ coordinates;

- if $\tilde{f}_t = \varphi_t^{-1} \circ f_t \circ \varphi_t$, then $\tilde{f}_t(r)$ is a local maximum of the $z$ coordinate restricted to $W^u(p_t)$ and also $\tilde{f}_t(r) = (t, q_0 + At)$;

- the least contracting direction of $\Lambda_t$ is the direction of $(1,0)$ in the linearizing coordinates restricted to $\{0\} \times \mathbb{R}^2$;

where $q_0 = (1,0)$ and $A \in \mathbb{R}^2$, see figure 3 where we have assumed $A = 0$.

We may extend the linearization neighborhood $L$ along $W^s(p_0)$ and $W^u(p_0)$ as explained in [9, Chap. 2]. In this way we may suppose that $r = f^{-N}q$ with $N = 1$. Moreover for small enough neighborhoods of $r$ and $q$ and for small $|t|$, $\tilde{f}_t$ has the form (see e.g. [10, Section 6])

$$(1 + z, X) \mapsto \left(t + az^2 + bX + h(t, z, X), \quad q_0 + At + Bz + CX + H(t, z, X)\right). \tag{8}$$

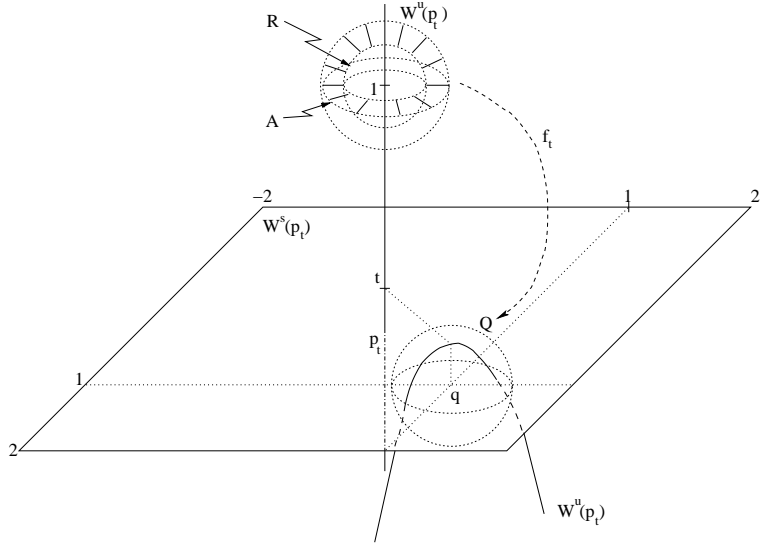

Figure 3: The adapted linearized coordinates and the regions of recurrence near the tangency points $q$ and $r$

Here $a \in \mathbb{R} \setminus \{0\}$, $A \in \mathbb{R}^2$, $B \in \mathbb{R}^2 \setminus \{0\}$, $b \in \mathcal{L}(\mathbb{R}^2, \mathbb{R}) \setminus \{0\}$, $C \in \mathcal{L}(\mathbb{R}^2, \mathbb{R}^2)$, $H$ is of order 2 or higher, and $h$ is of order 3 or higher in $z$ and of order 2 or higher in $t, X$ and $tz$. According to the assumptions of quadratic tangency and generic unfolding we must assume also that

$$Dh = 0, \quad DH = 0, \quad D_1 h = D_2 h = D_{12} h = 0 \quad \text{at} \quad (0, 0, 0^2). \tag{9}$$

Thus we may choose a compact neighborhood $\mathcal{Q}'$ of $q$ and $t_\star > 0$ small enough so that $\mathcal{Q}'$ together with the neighborhood

$$\mathcal{R} = \overline{\bigcup_{0 \leq t \leq t_\star} \tilde{f}_t^{-N} \mathcal{Q}'} \quad \text{of} \quad r = (1, 0, 0)$$

satisfy (8) for $\tilde{f}_t : \mathcal{R} \to \mathcal{Q}'$ and each $0 < t < t_\star$.

In addition, in the linearized domain it is easy to see that the action of $\tilde{f}_t$ on $\mathcal{Q} \subset \mathcal{Q}'$ will take every point $x \in \mathcal{Q}$ to $\mathcal{R}$ with an orbit inside $L$, for a compact neighborhood $\mathcal{Q}$ of $q$ inside $\mathcal{Q}'$ and for every $\underline{t} \in [0, t_\star]^\mathbb{N}$ — the number of iterates from $\mathcal{Q}$ to $\mathcal{R}$ will depend on $x \in \mathcal{Q}$ and $\underline{t}$.

We shall also take $\mathcal{R} \subset U$ where $U$ is a closed neighborhood of the basic set $\Lambda_p$ the saddle point $p$ belongs to. Moreover since $q \notin \Lambda_p$ we can make $\mathcal{Q}$ small enough so that $U \cap \mathcal{Q} = \emptyset$. For the arguments that follow we need to define a *small annulus* $\mathcal{A} = [\cup_{y \in \mathcal{R}} B(y, \zeta)] \setminus \mathcal{R}$, where $\zeta > 0$ is so small that $\mathcal{A} \subset U$ and $f_t(\mathcal{A}) \cap U = \emptyset$ for all $0 < t < t_\star$.

Finally, expression (7) guarantees that there is an invariant and expanding cone field $\mathcal{C}^u$ around the unstable direction $(1, 0, 0)$ defined on $T_x M$ at every $x \in L$. Moreover the stable foliation of $\Lambda_p$, being of codimension 1, can be extended to a $C^{1+\delta}$ foliation $\mathcal{W}$ of $U$, $0 < \delta < 1$, invariant by $f$ (see [6, Corollary 19.1.11] and [9, Appendix 1]) and we may assume that its leaves $\mathcal{W}(z, X)$ are "horizontal" hyperplanes $\{z = \text{constant}\}$ in the linearized coordinates. Since isolated hyperbolic sets admit *analytic continuations* [12] we have the same properties for the extension $\mathcal{W}_t$ of the stable foliation of $\Lambda_{p_t}$ to $U$ on linearized coordinates, where $p_t$ and $\Lambda_{p_t}$ are the hyperbolic continuations of $p$ and $\Lambda_p$ with respect to $f_t, 0 \leq t < t_\star$.

**Remark 3.1.** *The definition of $\mathcal{R}$ ensures that if $f_{\underline{t}}^k x \in \mathcal{Q}$ for some $k \geq 1$ and $x \in M$, then it must be that $f_{\underline{t}}^{k-1} x \in \mathcal{R}$.*

## 3.2 Returning to the neighborhood of the tangency

The linearization condition on the saddle fixed point $p$ and the condition of generic unfolding provides us with the following result.

**Lemma 3.2.** *There are $0 < b_0 < 1 < c_0$, small enough neighborhoods $\mathcal{Q}', \mathcal{Q}, \mathcal{R}$ and a small enough $t_\star > 0$ such that if $v \in T_x M$ with $x \in \mathcal{Q}$ and $t_1, \ldots, t_k \in [0, t_\star]$ is a sequence for which $k \in \mathbb{N}$ is the first return time of $x$ to $\mathcal{Q}$ (i.e. $f_{\underline{t}}^k x \in \mathcal{Q}$ and $f_{\underline{t}}^j x \notin \mathcal{Q}$ for $1 \leq j < k$), then*

$$\text{slope}(v) \geq c_0 \implies \begin{cases} \text{slope}(d(f_{\underline{t}}^k)_x v) \leq b_0, \\ \|d(f_{\underline{t}}^k)_x v\| \geq \frac{\|B\|}{10} \|v\|, \\ \angle(d(f_{\underline{t}}^k)_x v, B) \leq b_0. \end{cases}$$

Here and in the proof that follows the slope of a vector $u = (u_0, u_1, u_2) \in T_z M$ with $z \in L$ is to be understood has the ratio (after changing to the linearized coordinates) $|u_0|(\max\{|u_1|, |u_2|\})^{-1}$ and the norm $\|u\|$ equals $\max\{|u_0|, |u_1|, |u_2|\}$. The angle is calculated via the usual inner product in $\mathbb{R}^3$ in linearized coordinates. Hence this shows that vectors pointing to the expanding direction near $q$ will point to the contracting direction near $q$ after returning there under random iterations. Moreover these "return vectors" will point to a direction very close to that of $B$.

*Proof.* The previous adaptation of the linearized coordinates implies that in the stated conditions the sequence $f_{\underline{t}}, \ldots, f_{\underline{t}}^{k-1} x$ is in the linearization domain $L$ and also that $f_{\underline{t}}^{k-1} x \in \mathcal{R}$, where we have identified $f$ and the points in $L$ with the corresponding objects $\tilde{f}$ and points in the linear coordinates through $\varphi$. Then, assuming that $c_0$ is big enough, we must have that $v \in \mathcal{C}^u(x)$ and so $(df_{\underline{t}}^j)_x v \in C^u(f_{\underline{t}}^j x)$ for $j = 1, \ldots, k-1$. Moreover, by (7), the length will be expanded at an exponentially fast rate close to $\sigma$ and the slope will be incresed at a rate of $\eta = (\max\{\lambda_i \sigma, i = 1, 2\})^{-1} > 1$.

Assuming it takes at least $N$ iterations inside $L$ to go from $\mathcal{Q}$ to $\mathcal{R}$, we get that $w = (df_{\underline{t}}^{k-1})_x v$ satisfies $\|w\| \geq ((1+\sigma)/2)^N \|v\|$ and $\text{slope}(w) \geq \eta^N c_0$.

Let us now take $(1 + z, X) = f_{\underline{t}}^{k-1} x \in \mathcal{R}$. Writing everything in the linearized coordinates, we let $w = (w_0, W) \in \mathbb{R} \times \mathbb{R}^2$ and derive from (8) that the slope $[(df_{t_k})_z w]$

equals

$$\frac{\big|\,[2az + D_2h(t,z,X)] \cdot w_0 + [b + D_3h(t,z,X)] \cdot W\,\big|}{\max_{i=1,2}\{\big|\,[B_i + D_2H(t,z,X)] \cdot w_0 + [C_i + D_3H(t,z,X)] \cdot W\,\big|\}}$$

$$\leq \frac{\big|\,[2az + D_2h(t,z,X)]\,\big| + \big\|\,[b + D_3h(t,z,X)]\,\big\| \cdot [\text{slope}\,(w)]^{-1}}{\max_{i=1,2}\{\big|\,|[B_i + D_2H(t,z,X)]| - \|[C_i + D_3H(t,z,X)]\| \cdot [\text{slope}\,(w)]^{-1}\,\big|\}},$$

where $B_i \in \mathbb{R}$ and $C_i \in \mathcal{L}(\mathbb{R}^2, \mathbb{R})$ are the projections on the $i$th coordinate of $B \in \mathcal{L}(\mathbb{R}, \mathbb{R}^2)$ and $C \in \mathcal{L}(\mathbb{R}^2, \mathbb{R}^2)$.

Since $|2az|/|B_i|$ may be made smaller than any given positive number by making $\mathcal{R}$ smaller via shrinking $\mathcal{Q}', \mathcal{Q}$ and reducing $t_\star > 0$, we get that the slope above can be made smaller than any given $b_0 < 1$. Indeed, taking the neighborhoods $\mathcal{Q}'$ and $\mathcal{Q}$ smaller implies that the minimum number $N$ of iterates needed to take the points of $\mathcal{Q}$ to $\mathcal{R}$ grows without bound, so that the value of slope $(w)$ tends to infinity. Together with (9) this proves the statement for the slope.

We also have that $|w_0| = \|w\|$ since slope $(w) > 1$ and, moreover, there must be $i \in \{1,2\}$ such that $B_i \neq 0$. Hence the norm of the $i$th coordinate of $(df_{t_k})_z w$ is bounded from below by

$$\big|\,|[B_i + D_2H(t,z,X)] \cdot w_0| - |[C_i + D_3H(t,z,X)] \cdot W|\,\big| \geq \frac{|B_i|}{10} \cdot \|w_0\|,$$

because we can make slope $(w) \gg 1$. Since we may take $i$ such that $|B_i| = \|B\|$ (recall that this is the maximum norm), then $\|(df_{\underline{t}}^k)_x v\| \geq \|B\|\|w\|/10$ as long as $\mathcal{Q}', \mathcal{Q}, \mathcal{R}$ and $t_\star$ are small enough. Finally, if we calculate the usual inner product $(df_{t_k})_z w \cdot B$ and again use that slope $(w)$ can be made arbitrarily big, then we easily see that we can make $\angle((df_{t_k})_z w, B) \leq b_0$. $\square$

## 3.3 Iteration of the perturbation curve

Here we show that the perturbation at a return to $\mathcal{Q}$ acts essentially along the "vertical" unstable $z$ direction, due to the supposition of generic unfolding. Moreover assuming that every point of the curve traced by the set of all perturbations of a point has the same return time to $\mathcal{Q}$, we show that the return is also a curve but essentially parallel to the "horizontal" stable $X$ direction — this will be a consequence of the fold on the unstable manifold provided by the tangency.

Let $0 < b_0 < 1 < c_0$ be the constants obtained in the previous subsection.

**Lemma 3.3.** *Let $x \in M$, $\underline{t} \in T$ and $k \geq 1$ be such that $y = f_{\underline{t}}^{k-1} x \in \mathcal{R}$. Then the curve $\tilde{\gamma}(s) = f_s y, s \in \text{supp}\,(\theta)$ for some probability measure $\theta$ as in (2) satisfies*

$$\text{slope}\,(\tilde{\gamma}'(s)) \geq c_0 \quad \text{and} \quad \|\tilde{\gamma}'(s)\| \geq 1/2.$$

Let us assume that the return time of $\tilde{\gamma}(s)$ to $\mathcal{Q}$ under a perturbation vector $\underline{v}$ is $R > 1$ for all $s \in supp\,(\theta)$ and $\underline{v} \in supp\,(\theta^{\mathbb{N}})$. Then the map $\gamma(u,s) = f_u(f_{\underline{v}}^{R-1}\tilde{\gamma}(s))$, $s, u \in supp\,(\theta)$ is a diffeomorphism onto its image satisfying for all $s, u \in supp\,(\theta)$:

$$\text{slope}\,(\partial_u \gamma) \geq c_0, \qquad \|\partial_u \gamma\| \geq \frac{1}{2}; \tag{10}$$

$$\text{slope}\,(\partial_s \gamma) \leq b_0, \quad \text{and} \quad \|\partial_s \gamma\| \geq \frac{\|B\|}{10}. \tag{11}$$

*Proof.* The result on the slope and the norm of $\tilde{\gamma}'(s)$ follow from the adaptation of the linearizing coordinates done in subsection 3.1, as long as $t_\star > 0$ is taken small enough.

The assumption of constant return time $R$ for all points of the trace of $\tilde{\gamma}$ implies that Lemma 3.2 can be applied to each $\tilde{\gamma}(s)$, providing the bounds on the slope and norm of $\partial_s \gamma$. Since $\gamma(u,s)$ for fixed $s$ is in the conditions of $\tilde{\gamma}$, we get also the bounds on the slope and norm of $\partial_u \gamma$ as before. □

**Remark 3.4.** *Taking $\mathcal{Q}$ and $t_\star$ small enough, we can ensure that $\partial_s \gamma$ is very close to the direction of $B$, the tangent direction of $W^u(p)$ at $q$, i.e. the angle is bounded $\angle(B, \partial_s \gamma) \leq b_0$ by Lemma 3.2. Moreover under the same conditions $\partial_u \gamma$ is very close to the direction of the $z$ axis in the linearizing coordinates because of its adaptation (8) and (9). In addition, since the speed of the curves $\gamma(\cdot, s)$ is bounded from below, we see that there is a nearly vertical curve in the image of all perturbations at a return time with length proportional to the diameter of $supp\,(\theta)$.*

**Remark 3.5.** *Under the assumptions of this lemma, since the bounds obtained are uniform we deduce that the set of all perturbed iterates of $x$ up to time $k + R$, $\{f_{\underline{t}}^{k+R} x : \underline{t} \in supp\,(\theta^{\mathbb{N}})\}$, contains a 2-disk given by the image of $\gamma$ above. Note that this 2-disk is formed by a family of nearly vertical curves $\gamma_s(u) = \gamma(u, s)$. Moreover the bounds provided by Remark 3.4 ensure that the diameter of this disk will be greater than a constant times the diameter of $supp\,(\theta)$, that is, of the form $K\epsilon$ independent of $x$ and the number of iterates $k + R$ involved in the construction, as long as $t_\star > 0$ is small.*

## 4 Finite number of physical measures

In what follows we fix a small enough $t^* > 0$, a point $t_0 \in ]0, t^\star[$ and a noise level $\epsilon \in ]0, t^\star - t_0[$ as in the statement of Theorem A.

### 4.1 Regular orbits near the tangency

In this subsection we prove a result about $\omega$ limit points of generic random orbits. Recall that $(\theta_\epsilon)_{\epsilon>0}$ is a family of probability meausures given by (2) (see also Remark 1.6).

**Definition 4.1.** *Given $x \in M$ and $\underline{t} \in T$ we let $\omega(x, \underline{t})$ be the closure of the set of accumulation points of the sequence $(f_{\underline{t}}^n x)_{n \geq 1}$.*

**Proposition 4.2.** *Let $x$ be a $\theta_\epsilon$ recurrent point with respect to a set $W \subset \text{supp}(\theta_\epsilon^\mathbb{N})$. Then for $\theta_\epsilon^\mathbb{N}$ almost all $\underline{t} \in W$ every $z \in \omega(x, \underline{t})$ is a regular point, i.e. $r(k, \underline{s}, z) = r(k, z) \leq J$ for all $k \geq 2$ where $J \in \mathbb{N}$ is independent of $x, W$ and $z$, and the return times do not depend on $\underline{s} \in \text{supp}(\theta_\epsilon)$.*

*Proof.* Let us take $y \in \omega(x, \underline{t})$ for $x$ a $\theta = \theta_\epsilon$ recurrent point with respect to $W \subset T$, where $\theta^\mathbb{N}(W) > 0$ and $\underline{t} \in W$ is $\theta^\mathbb{N}$ generic. Now we proceed as follows.

- If $f_{\underline{s}}^k y \in \mathcal{Q}$ for some $\underline{s} \in \text{supp}(\theta) \in T$ and $k \geq 1$, then $f_{\underline{u}}^k y \in \mathcal{Q}$ for all $\underline{u} \in \text{supp}(\theta)$.

In fact, if $f_{\underline{u}}^k y \notin \mathcal{Q}$ for some $\underline{u} \in \text{supp}(\theta)$, then the construction of $\mathcal{Q}$ and $\mathcal{R}$ forces $f_{\underline{s}}^{k-1} y \in \mathcal{R}$ and $f_{\underline{u}}^{k-1} y \notin \mathcal{R}$, cf. Remark 3.1. Since $\text{supp}(\theta)$ is connected and $T \ni \underline{v} \mapsto f_{\underline{v}}^{k-1} y$ is continuous, we know that there is $\underline{v} \in \text{supp}(\theta)$ such that $f_{\underline{v}}^{k-1} y \in \mathcal{A}$. This implies that $w = f_{\underline{v}}^k y \in M \setminus (U \cup \mathcal{Q})$.

However, since $y \in \omega(x, \underline{t})$ with $\underline{t}$ a $\theta_\epsilon^\mathbb{N}$ generic point, we may get for small $\delta > 0$ a $\theta^\mathbb{N}$ generic vector $\underline{\omega} \in W$ and $n$ big enough such that:

1. $d(\underline{\omega}, \underline{t}) < \delta$ and $d(\sigma^n \underline{\omega}, \underline{v}) \leq \delta$ by Lemma 2.4;

2. $\text{dist}(f_{\underline{\omega}}^n x, y) < \delta$ using the definition of $y$;

(where $d$ is the usual distance in $T$, $d(\underline{s}, \underline{t}) = \sum_{j \geq 1} 2^{-j} |s_j - t_j|$, and dist is the Riemannian induced distance on $M$). Now because $\mathcal{A}$ is open, if we take $\delta$ small enough, then combining the distances above we arrive at $f_{\underline{\omega}}^{n+k-1} x \in \mathcal{A}$ and thus $f_{\underline{\omega}}^{n+k} x \in M \setminus (U \cup \mathcal{Q})$. This contradicts the assumed $\theta$ recurrence of $x$ with respect to $W$.

- $f_{\underline{s}}^k y \in U \cup \mathcal{Q}$ for all $k \geq 1$.

This is obvious from the arguments proving the previous item above.

- If $f_{\underline{s}}^k y \in U$ for all $k \geq k_0$ for some $k_0 \in \mathbb{N}$ and $\underline{s} \in \text{supp}(\theta)$, then $f_{\underline{u}}^k y \in U$ for all $k \geq k_0$ and every $\underline{u} \in \text{supp}(\theta)$.

Indeed, let us assume that there exists $\underline{u} \in \text{supp}(\theta)$ and $l \geq k_0$ such that $f_{\underline{u}}^l y \notin U$. Since $f_{\underline{s}}^l y \in U$, the connectedness argument above likewise implies that there is $\underline{v} \in \text{supp}(\theta)$ for which $f_{\underline{v}}^l y \in M \setminus (U \cup \mathcal{Q})$. But now we can repeat the arguments from the previous items.

After these items we know that the orbit of $y$ under any perturbation $\underline{s} \in \text{supp}(\theta)$ falls into one of the following cathegories:

1. there is a finite number of passages through $\mathcal{Q}$ and afterwards wanders in $U$;

2. always remains in $U$ and never passes through $\mathcal{Q}$;

3. passes through $\mathcal{Q}$ infinitely many times and the passage times do not depend on $\underline{s}$.

We now eliminate the first two possibilities.

- *Finite number of returns only is not possible.*

Assume that $r \in \mathbb{N}$ is the last return iterate of $y$ to $\mathcal{Q}$ under every $\underline{s} \in \text{supp}(\theta)$ — the passage times of the orbits through $\mathcal{Q}$ do not depend on $\underline{s}$, according to the preceeding arguments in this section. Then we know from Lemma 3.3 and Remark 3.4 that the set $\{f^r_{\underline{s}} y : \underline{s} \in \text{supp}(\theta)\}$ contains a curve $c : (-\epsilon, \epsilon) \to \mathcal{Q}$ with slope $\geq c_0$ and speed $\geq \frac{1}{2}$. In the future, the points $c(t)$ of the curve will stay forever in $U$ under every perturbation, i.e. $f^k_{\underline{u}}(c(t)) \in U$ for all $\underline{u} \in \text{supp}(\theta), k \geq 1$ and $t \in (-\epsilon, \epsilon)$. But $U$ is an isolating neighborhood of a hyperbolic set and $c$ is a curve inside the unstable cone field, thus the curves $c_k = f^k_{\underline{s}} \circ c$ cannot remain inside $U$ for all $k \geq 1$ because their lenght will tend to infinity while staying always inside the unstable cone field.

- *Never returning is not possible.*

Since we are considering $y = (z, X) \in \omega(x, \underline{t})$ with a $\theta$ recurrent $x$ we may assume that $y \in \mathcal{Q}$. If $f^k_{\underline{s}} y \in U$ for all $k \geq 1$ and $\underline{s} \in \text{supp}(\theta)$, then by a similar reazoning to that of the previous subsection it is not possible that, for any given $\delta > 0$, the segment $\gamma = \{(z + r, X) : -\delta < r < \delta\}$ satisfies $f^k_{\underline{s}} \gamma \subset U$ for all $k \geq 1$ and any $\underline{s} \in \text{supp}(\theta)$. Thus there must be two open subsegments (since $U \cup \mathcal{Q}$ is closed) of $\gamma$, one $\gamma^+$ "above" and the other $\gamma^-$ "below" $y$ such that there are $k^+, k^- \in \mathbb{N}$ satisfying (recall that $t_0 \in \text{supp}(\theta)$)

$$(f^{k^+}_{t_0} \gamma^+ \cup f^{k^-}_{t_0} \gamma^-) \subset M \setminus (U \cup \mathcal{Q}).$$

We will assume that $k^\pm$ were taken minimal with respect to this property. Now for $w = (z + r, X) \in \gamma^\pm$ we may consider the leaf $\mathcal{W}_{t_0}(w)$ through $w$. Hence $\text{dist}(f^k_{t_0} w, f^k_{t_0} w_1) \leq C\lambda^k \|w - w_1\|$ for $1 \leq k < k^\pm$, where $\lambda \in (\lambda_1, 1)$ and $C > 0$ is a constant bounding the distortion induced by the linearizing coodinate change. Noting that we can make $k_1 = \max\{k^+, k^-\}$ bigger by taking $\delta$ smaller, we see that there are two "plates" $D^\pm$ made by leaves through points in $\gamma^\pm$ having a radius $\leq C^{-1} \lambda^{-k_1} d_0$, where $2d_0 = \text{dist}(U, \mathcal{Q}) > 0$, see figure 4. Points in any of these plates will leave $U$ after no more than $k_1$ iterates under $f_{t_0}$.

There are $n_1 < n_2 < \ldots$ such that $f^{n_k}_{\underline{t}} x \to y$ when $k \to +\infty$ by definition of $y$. We remark that under our assumptions this implies $f^{n_k}_{\underline{t}} x \in \mathcal{Q}$. We know from Lemma 3.3 that there is a curve $c : \text{supp}(\theta) \to \mathcal{Q}$ given by $c(u) = f_u(f^{n_k - 1}_{\underline{t}} x)$ whose slope at every point is $\geq c_0$ and whose speed is bounded away from zero. Hence by taking $\delta$ smaller and $k$ big enough we arrive at the situation depicted in figure 4. We notice that one of $D^\pm$ intersects the trace of $c$ in a segment $I$ of length $\geq l_0$, because the slope of $c$ is bounded from below by a constant $c_0 \gg 1$. In addition, since the speed of $c$ is also bounded away from zero by a uniform constant, then $I_0 = c^{-1}(I) \subset \text{supp}(\theta)$ is a segment of length $\geq l_1$ for some $l_1 > 0$ depending neither on $x$ nor on $n_k$. Hence we have either

$$f^{k^+}_{t_0}\big(c(I_0)\big) \subset M \setminus (U \cup \mathcal{Q}) \quad \text{or} \quad f^{k^-}_{t_0}\big(c(I_0)\big) \subset M \setminus (U \cup \mathcal{Q}).$$

However for every $\gamma, \delta \in (0, 1)$ there exists $k_0 \in \mathbb{N}$ such that for $k \geq k_0$ there are subsets $W_{n_k - 1} \subset W$ and $Y_{n_k - 1} \subset \text{supp}(\theta)$ satisfying items 1 through 4 of Lemma 2.7. Since $\theta(I_0)$ is bounded away from zero independently of $k$, we see that we can find $\underline{s} \in Y_{n_k-1}$ very close to the vector $(t_0, t_0, \dots)$ (taking $\gamma$ close to 1) such that $\theta\big(I_0 \cap p_{n_k} W_{n_k-1}(\underline{t}, n_k-1, \underline{s})\big) > 0$, as long as $\delta$ is close enough to 1. Thus there is in $W$ a subset of positive $\theta$ measure whose vectors send $x$ into $M \setminus (U \cup \mathcal{Q})$, which contradicts the standing assumption on $x$ and $W$.

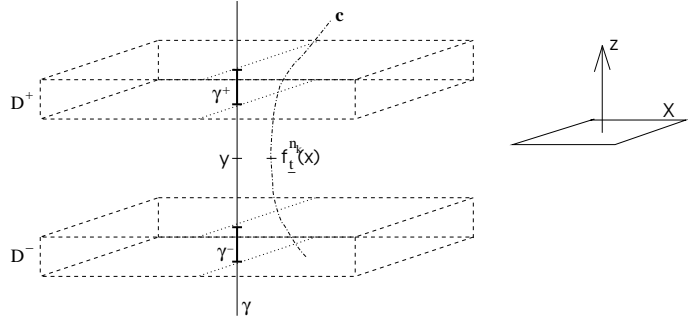

Figure 4: The stable plates and the curve $c$ near $y$.

Hence the orbit $(f^n_{\underline{s}} z)_{n \geq 1}$ has infinitely many returns to $\mathcal{Q}$ at iterates that do not depend on $\underline{s} \in \operatorname{supp}(\theta)$.

- *The uniform bound on the return times.*

We first wait for the first return of $(f^n_{\underline{s}} z)_{n \geq 1}$ to $\mathcal{Q}$ for $\underline{t} = (t_0, t_0, \dots)$ say. Note that by doing this we loose control on the number of iterates until the first return — this may depend on $z$.

At the return, the perturbation will "generate" a nearly "vertical" curve $c$ in $\mathcal{Q}$, whose length is bounded from below by $K\epsilon$, as consequence of Lemmas 3.2 and 3.3, where $K > 0$ is a uniform constant. Since every point of $c$ will again return to $\mathcal{Q}$ after the same number of iterates $n \geq 1$, it must be that $f^n_{\underline{t}}(c) \subset \mathcal{R}$. But $c, f_{\underline{t}}(c), \dots, f^n_{\underline{t}}(c) \subset U$ and so the curve $c$ will be uniformly streched in the expanding direction during these $k$ iterations, thus $K\epsilon\sigma^n \leq \operatorname{diam} \mathcal{R}$.

This shows that there must be a uniform constant $J \in \mathbb{N}$ bounding $n$ from above. Hence we arrive at $r(k, \underline{t}, z) \leq J$ for all $k \geq 2$. □

## 4.2 The basin of ergodic measures near the tangency

Now we show that the 2-disk $\Delta$ obtained in subsection 3.3 returns to $\mathcal{Q}$ as another 2-disk, but whose direction is transversal to $\Delta$ due to the generic condition 2(d) on the quadratic tangency. The perturbation at this return will then generate a nearly vertical curve through each point of the return disk, providing a ball in $\mathbb{R}^3$ — see figure 5 — and also the absolute continuity of the push-forwards (5).

To arrive at this result the return times must be independent of the perturbation chosen, so we shall consider regular orbits in the arguments that follow. Moreover we will use linearized coordinates throughout.

### 4.2.1 Vertical 2-disks and their returns

We know from subsection 3.3 that if we take $x \in M$ a $\theta_\epsilon$ regular point as given in Proposition 4.2, then the set $\{f_{\underline{t}}^k x : \underline{t} \in \text{supp}\,(\theta_\epsilon^{\mathbb{N}})\}$ contains a nearly vertical 2-disk $\Delta$, where $k = r(i,x)$ is the $i$th return time to $\mathcal{Q}$ with $i \geq 2$. The tangent space to each point of $\Delta$ is generated by a vector very close to the expanding "vertical" direction and by a vector very close to the direction of the vector $B$ — see Remarks 3.4 and 3.5.

The points of $\Delta$ are all $\theta_\epsilon$ regular, since they are images of a $\theta_\epsilon$ regular point under a random requence of maps $f_{\underline{t}}^k$. Hence we can consider the return of each point of $\Delta$ to $\mathcal{Q}$ under the action of $f_{t_0}$. Let us assume that the number of iterates until the return is $R \geq 1$. Then $\Delta' = f_{t_0}^{R-1}(\Delta) \subset \mathcal{R}$ because $\Delta$ consists of regular points.

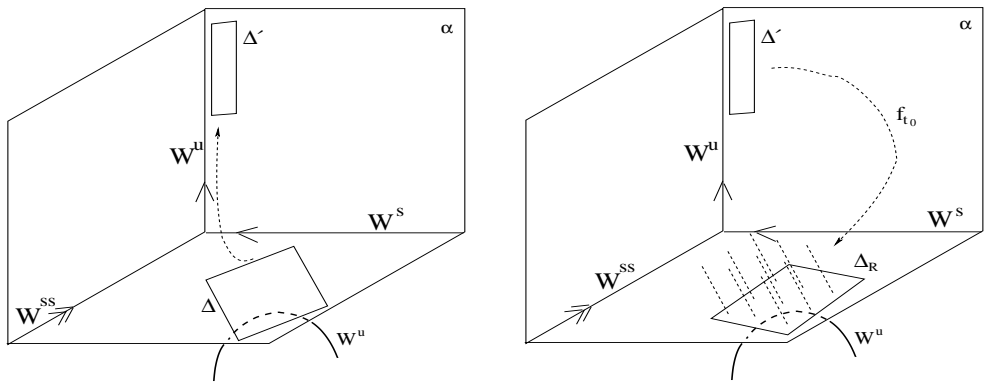

Figure 5: Sketch of the directions of the disks $\Delta$, $\Delta'$ and $\Delta_R$ relative to the invariant directions $W^{ss}, W^s, W^c$ and the perturbation curves through $\Delta_R$.

Before arriving at $\mathcal{R}$, the disk $\Delta$ was expanded and contracted in a finite number of hyperbolic iterates. Then $\Delta'$ also is a nearly vertical 2-disk since the expanding direction is the "vertical" one.

However $\Delta'$ is also a disk "close" to the direction of the plane $\alpha$ defined by the expanding $z$ axis and the weak contracting direction, where we interpret *planes with close directions* in the usual intuitive way of the angle between the directions of their orthogonal lines.

Indeed, since the tangent direction to the unstable manifold at the tangency (the direction of $B$) is transversal to the strong-stable direction of the saddle fixed point $p_0$, the persistence of this configuration under perturbation ensures that the components of any tangent vector $v$ to $\Delta$, in the strong-stable direction, will suffer a much greater contraction than the components of the same vector in the weak-stable direction — see the left hand side of figure 5.

Thus we can guarantee that the tangent directions to every point of the disk $\Delta_R = f_{t_0}(\Delta')$ will be close to the tangent directions to the $f_{t_0}$-image of the projection of $\Delta'$ on the plane $\alpha$, in the strong-stable direction — see the right hand side of figure 5. Moreover any degree of closeness can be achieved (i.e. sufficiently small angles) by shrinking $\mathcal{Q}$ and $\mathcal{R}$ without affecting the constants in the arguments of Lemma 3.3, Remarks 3.4 and 3.5.

Finally the generic condition 2(d) on the family $\Phi$ clearly implies that $\Delta_R$ is transversal to the unfolding direction $D$ near $q$. Since $\Delta_R$ is obtained after iteration under $f_{t_0}$, we can now consider the set of all perturbations and use Lemma 3.3 to see that this set contains a nearly "vertical" curve through each point of $\Delta_R$ with length $\geq K\epsilon$, cf. figure 5, Lemma 3.3 and Remarks 3.4 and 3.5. Thus there must be a ball of uniform radius (dependent on $\Phi$ and $\epsilon$ only) inside the set of all perturbed iterates of a regular point after three consecutive returns.

We note that since for regular points the number of iterates between consecutive returns to $\mathcal{Q}$ is bounded from above, then in fact there is a uniform bound for the radius of the inner ball for all iterates beyond $r(1, x) + r(2, x) + r(3, x)$ (the third return time of $x$ to $\mathcal{Q}$).

### 4.2.2 Absolute continuity of the push forward

In addition, let $R(k, x) = \sum_{l=1}^{k} r(l, x)$ be the $k$th return iterate and $i = R(1, x) < j = R(2, x) < k = R(3, x)$ be three successive return times to $\mathcal{Q}$. Then the map $\phi(t_1, \ldots, t_k) = (f_{t_k} \circ \cdots \circ f_{t_1})x$, for $t_1, \ldots, t_k \in \mathrm{supp}\,(\theta_\epsilon)$, is a submersion because $D_{i,j,k}\phi : \mathbb{R}^3 \to \mathbb{R}^3$ is a isomorphism.

This forces $(f_x^n)_*\theta_\epsilon^\mathbb{N} \ll m$ for all $n \geq k$. Indeed, we have just obtained this for $n = k$ since $(f_x^k)_*\theta_\epsilon^\mathbb{N} = \phi_*\theta_\epsilon^\mathbb{N}$. For higher $n$ we note that if $E \subset M$ with $m(E) = 0$, then $[(f_x^n)_*\theta_\epsilon^\mathbb{N}](E) = \int \theta_\epsilon^k(E(t_{k+1}, \ldots, t_n))\, d\theta_\epsilon^{n-k}(t_{k+1}, \ldots, t_n)$ by Fubini's Theorem, where $E(t_{k+1}, \ldots, t_n) = \{(t_1, \ldots, t_k) : \phi(t_1, \ldots, t_k) \in (f_{t_n} \circ \cdots \circ f_{t_{k+1}})^{-1}E\}$. But since each $f_t$ is a diffeomorphism, we get $m((f_{t_n} \circ \cdots \circ f_{t_{k+1}})^{-1}E) = 0$ and so $\theta_\epsilon^k(E(t_{k+1}, \ldots, t_n)) = [(f_x^k)_*\theta_\epsilon^\mathbb{N}](f_{t_n} \circ \cdots \circ f_{t_{k+1}})^{-1}E = 0$.

Hence we have proved

**Proposition 4.3.** *Let $x$ be a regular point for the noise level $\epsilon$. Then the set $\{f_{\underline{t}}^n x : \underline{t} \in \mathrm{supp}\,(\theta_\epsilon^\mathbb{N})\}$ contains a ball with radius $\geq K\epsilon$ centered around $f_{t_0}^n x$, where $K$ depends on $\Phi$ and $\epsilon$ only, for every $n \geq R(3, x)$. Moreover it holds also that $(f_x^n)_*\theta_\epsilon^\mathbb{N} \ll m$ for all $n \geq R(3, x)$.*

## 4.3 Nonempty interior and ergodic decomposition

Let $\mu$ be a stationary ergodic probability measure whose support is inside $U \cup \mathcal{Q}$ and such that $\mu(\mathcal{Q}) > 0$. Since $\theta_\epsilon^\mathbb{N} \times \mu$ is $F$ invariant (see section 2.1), Poincaré's Recurrence Theorem ([8, Chpt. 1.2]) ensures that $\theta_\epsilon^\mathbb{N} \times \mu$ almost all $(\underline{t}, x) \in \mathrm{supp}\,(\theta_\epsilon^\mathbb{N} \times \mu)$ is $\omega$ recurrent with respect to the action of $F$. But then $x$ is a $\theta_\epsilon$ recurrent point with respect to $W \subset T$ with $\theta_\epsilon^\mathbb{N}(W) = 1$ and, moreover, $x \in \omega(x, \underline{t})$ for all $\underline{t} \in W$.

Hence there exists a regular point $x \in \mathrm{supp}\,(\mu)$ and so Proposition 4.3 shows that the interior of $\mathrm{supp}\,(\mu)$ is nonempty and contains a ball of uniform radius.

Let $\mu_1, \mu_2$ be two distinct stationary ergodic probability measures. Since ergodic basins are always disjoint it holds that $B(\theta_\epsilon^\mathbb{N} \times \mu_1) \cap B(\theta_\epsilon^\mathbb{N} \times \mu_2) = \emptyset$ and because $\mu_i(B(\theta_\epsilon^\mathbb{N} \times \mu_i) \setminus \mathrm{supp}\,(\theta_\epsilon^\mathbb{N} \times \mu_i)) = 0$ by the Ergodic Theorem, we deduce that $\mu_i(\mathrm{supp}\,(\mu_j)) = \delta_{ij}$, $i, j = 1, 2$.

This means that $\mu_1 \perp \mu_2$. So the number of stationary ergodic probability measures supported in $U \cup \mathcal{Q}$ and giving positive weight to $\mathcal{Q}$ is finite, because each support of one of these measures contains a different ball of uniform radius and volume unifomly bounded from below, and $m(U \cup \mathcal{Q}) < 1$. This proves item 1 of Theorem A.

Together with the ergodic decomposition result for stationary probability measures in section 2.1 we have proved

**Proposition 4.4.** *For every stationary probability measure $\mu$ there are $\alpha_1, \ldots, \alpha_l \geq 0$ such that $\alpha_1 + \cdots + \alpha_l = 1$ and $\mu = \sum_{i=1}^{l} \alpha_i \mu_i$.*

## 4.4 Absolute continuity

Let us consider the connected components of the interior $\mathcal{U}$ of $\operatorname{supp}(\mu)$ for a given stationary ergodic probability measure $\mu$. Since $\mathcal{U}$ is invariant in the sense of (4), these components must be permuted by the action of $\Phi(t)$ independently of $t \in \operatorname{supp}(\theta_\epsilon)$, because $\Phi(t)$ varies continuously with $t$ and $\operatorname{supp}(\theta_\epsilon)$ is also assumed to be connected (cf. Remark 1.6). Moreover the permutation is cyclic because $\mu$ is ergodic and, since these components are open sets inside $\operatorname{supp}(\mu)$, there is a $\omega$ recurrent and regular point inside each component by the arguments of the previous subsection.

Using again Proposition 4.3 we deduce that each connected component of $\mathcal{U}$ must contain a ball of uniform radius. Thus there are only finitely many connected components of $\mathcal{U}$ and the points in these components must have orbits whose return times to $\mathcal{Q}$ are uniformly bounded from above.

Hence in fact *every point $x \in \mathcal{U}$ is regular*! This means that it exists $J \in \mathbb{N}$ such that there are *three returns* to $\mathcal{Q}$ of the orbit of all $x \in \mathcal{U}$ in fewer than $J$ iterates.

We can consider the normalized restriction $\mu_\mathcal{U}$ of $\mu$ to $\mathcal{U}$. This is a stationary probability measure by [2, Lemma 8.2] and the ergodicity of $\mu$ implies that $\mu = \mu_\mathcal{U}$. Thus $\mu(\partial \operatorname{supp}(\mu)) = 0$ and all the mass of $\mu$ is in $\mathcal{U}$.

Now iterating equality (3) we easily obtain $\mu\varphi = \int [(f_x^J)_* \theta_\epsilon^\mathbb{N}] \varphi \, d\mu(x)$ for any $\varphi \in C^0(M, \mathbb{R})$. Since $C^0(M, \mathbb{R})$ is dense in $L^1(\mu)$ in the $L^1$ norm, this identity is valid for all $\mu$ integrable real functions. Finally take $\varphi = 1_E$ where $E \subset M$ with $m(E) = 0$ and use the absolute continuity from Proposition 4.3 to conclude that $\mu(E) = 0$.

We have obtained

**Proposition 4.5.** *Every stationary ergodic probability measure $\mu$ with $\operatorname{supp}(\mu) \subset U \cup \mathcal{Q}$ and $\mu(\mathcal{Q}) > 0$ is absolutely continuous (with respect to $m$). Moreover the mass of $\mu$ lies in $\mathcal{U} = \operatorname{int}(\operatorname{supp}\mu)$ and $\mathcal{U}$ is the union of finitely many connected components cyclically permuted by the action of $\Phi(t)$, for all $t \in \operatorname{supp}(\theta_\epsilon)$.*

*In particular the return times to $\mathcal{Q}$ for every $x \in \mathcal{U}$ are uniformly bounded and independent of the perturbation.*

This proves item 2 of Theorem A and from the absolute continuity we obtain also that these ergodic probability measures are physical, by the Ergodic Theorem.

## 4.5 The orbits of recurrent points

Here we prove the existence of the partition of $supp\,(\theta_\epsilon^\mathbb{N})$ for all $\theta_\epsilon$ recurrent points stated in item 3 of Theorem A. We start with two preparatory lemmas.

The setting for next result is exactly provided by Propositions 4.3 and 4.5.

**Lemma 4.6.** *The ergodic basin $B(\mu)$ of every stationary ergodic probability measure $\mu$ with $supp\,(\mu) \subset U \cup \mathcal{Q}$ and $\mu(\mathcal{Q}) > 0$ contains an open set: $B(\mu) \supset \mathcal{U}$ where $\mathcal{U} = int\,(supp\,\mu)$.*

*Proof.* See the proof of Lemma 7.2 in [2]. □

This in turn gives the following.

**Lemma 4.7.** *For $x \in M$ let $W_i(x) = \{\underline{t} \in supp\,(\theta_\epsilon^\mathbb{N}) : n_k^{-1} \sum_{j=0}^{n_k-1} \delta_{f_{\underline{t}}^j x} \xrightarrow{w^\star} \mu_i\}$ for some given integer sequence $1 \leq n_1 < n_2 < \ldots$ and let $W_{i,j}(x) = \{\underline{t} \in supp\,(\theta_\epsilon^\mathbb{N}) : f_{\underline{t}}^j x \in int\,B(\mu_i)\}$, where $j \geq 1$ and $i = 1, \ldots, l$. Then $W_i(x) = \cup_{j\geq 1} W_{i,j}(x)$ modulo a $\theta_\epsilon^\mathbb{N}$ null set.*

**Remark 4.8.** *This shows that $W_i(x)$ does not depend on the given sequence $(n_k)_{k\geq 1}$.*

*Proof.* On the one hand, let us fix $i \in \{1, \ldots, l\}$, $\underline{t} \in W_i(x)$ and $\mathcal{U} = int\,(supp\,\mu_i)$. We take $\varphi \in C^0(M, \mathbb{R})$ with $\varphi \geq 0$ and $supp\,(\varphi) \subset V_y \subset B(\mu_i)$, where $V_y$ is an open neighborhood of a point $y \in \mathcal{U}$. Since $\mu_i \varphi > 0$, by the definition of $W_i(x)$ and the form of $\varphi$ there must be $j \geq 1$ such that $f_{\underline{t}}^j x \in V_y$. Hence $W_i(x) \subset \cup_{j\geq 1} W_{i,j}(x)$.

On the other hand, let $\underline{t} \in W_{i,j}(x)$ for some fixed $j \geq 1$. Since $\underline{t} \mapsto f_{\underline{t}}^j x$ is continuous, there is a neighborhood $W$ of $\underline{t}$ in $supp\,(\theta_\epsilon^\mathbb{N})$ such that $f_{\underline{s}}^j x \in int\,B(\mu)$ for all $\underline{s} \in W$. Hence $\theta_\epsilon^\mathbb{N}$ almost all $\underline{s} \in W$ are in $W_i(x)$. □

Now we are ready to obtain the partition. Let $x \in M$ be a $\theta_\epsilon$ recurrent point with respect to $W \subset supp\,(\theta_\epsilon^\mathbb{N})$. Then $y \in \omega(x, \underline{t})$ is a regular point for $\theta_\epsilon^\mathbb{N}$ almost every $\underline{t} \in W$ by Proposition 4.2. Thus every weak* accumulation point $\mu$ of the sequence $(\mu_n(y))_{n\geq 1}$ is a stationary probability measure (recall Lemma 2.1) satisfying $supp\,(\mu) \subset U \cup \mathcal{Q}$ and $\mu(\mathcal{Q}) \geq J^{-1}$ by definition of regular point. Hence by Proposition 4.4 we may decompose $\mu = \sum_{i=1}^l \alpha_i \mu_i$ with $\alpha_i \geq 0$ and $\sum_{i=1}^l \alpha_i = 1$.

**Lemma 4.9.** *If $\alpha_i > 0$ for some $i \in \{1, \ldots, l\}$ then both $\theta_\epsilon^\mathbb{N}(W_i(y)) > 0$ and $\theta_\epsilon^\mathbb{N}(W_i(x)) > 0$.*

*Proof.* Just take $\varphi \in C^0(M, \mathbb{R})$ with $\varphi \geq 0$ and $supp\,(\varphi) \subset V_z \subset B(\mu_i)$, where $V_z$ is an open neighborhood of a point $z \in int\,(supp\,\mu_i)$. Then $\mu\varphi \geq \alpha_i \mu_i(\varphi) > 0$ by construction of $\mu$ and so there must be $j \geq 1$ such that $\theta_\epsilon^\mathbb{N}\big((f_y^j)^{-1} V_z\big) > 0$. This implies that $\theta_\epsilon^\mathbb{N}(W_{i,j}(y)) > 0$ and $\theta_\epsilon^\mathbb{N}(W_i(y)) > 0$ follows by Lemma 4.7.

Since $y \in \omega(x, \underline{t})$ with $\underline{t} \in W$ a $\theta_\epsilon^\mathbb{N}$ generic vector, we take $V = (f_y^j)^{-1} V_z$ and use Lemma 2.4 to obtain:

- $k_0 \in \mathbb{N}$ such that $\theta_\epsilon^\mathbb{N}(\sigma^{-k} V \cap B_\eta(\underline{t}) \cap W) > 0$ for every $k \geq k_0$ and $\eta > 0$ (here the distance in $T$ is defined in Section 2.2);

- a sufficiently high $k \geq k_0$ and a sufficiently small $\eta > 0$ such that $f_{\underline{t}}^k x$ is close enough to $y$, so that $f_{\underline{s}}^{j+k} x \in V_z$ for all $\underline{s} \in (\sigma^{-k} V \cap B_\eta(\underline{t}) \cap W)$.

Then we get $\theta_\epsilon^{\mathbb{N}}\big((f_x^{j+k})^{-1} V_z\big) > 0$ and from this $\theta_\epsilon^{\mathbb{N}}(W_i(x)) > 0$ follows easily. $\square$

The sets $W_1(x), \ldots, W_l(x)$ form a partition of $W$ modulo $\theta_\epsilon^{\mathbb{N}}$ null sets since they are pairwise disjoint and their union equals $W$. Indeed, if we suppose that $W_0 = W \setminus (W_1(x) \cup \cdots \cup W_l(x))$ is such that $\theta_\epsilon^{\mathbb{N}}(W_0) > 0$, then $x$ is a $\theta_\epsilon$ recurrent point with respect to $W_0$. Repeating the arguments leading to the construction of $W_1(x), \ldots, W_l(x)$, we obtain $\theta_\epsilon^{\mathbb{N}}(W_0 \cap W_i(x)) > 0$ for some $i = 1, \ldots, l$, contradicting the definition of $W_0$. This concludes the proof of Theorem A.

**Remark 4.10.** *Is not difficult to see that the coefficients $\alpha_1, \ldots, \alpha_l$ of the decomposition of a stationary probability measure $\mu$, given as a weak\* accumulation point of $(\mu_n(x))_n$, satisfy $\alpha_i = \theta_\epsilon^{\mathbb{N}}(W_i(x))$. Hence the weak\* limit $\lim_n \mu_n(x)$ exists, see [3, Lemma 3.4].*

# 5 Homoclinic tangencies in higher dimensions

In the setting of manifolds with dimension $\geq 4$, we can repeat all the arguments in Sections 2, 3 and 4.1 in an obvious fashion, up to the point where we obtain three independent directions after three consecutive returns to $\mathcal{Q}$, in Section 4.2.

Form here on it is not clear how to obtain another independent direction at a return time to $\mathcal{Q}$. Indeed, the existence of a least contracting direction forces the three independent directions obtained up until now to approach the bidimensional plane given by the unstable and weakly constracting directions — as in Figure 5. At the next return the set of previous return images will be very close the disk $\Delta_R$.

We observe that the least contracting direction was strongly used in [10] and [11] to allow some bidimensional arguments to be applied through a reduction of codimension.

Hence though we might need to eliminate condition 1(c) of existence of a least contracting direction, doing so we may loose the main motivation behind this problem: the coexistence of infinitely many attractors, because the proofs in [10] possibly no longer hold. Moreover we can also loose the ability to apply some bidimensional arguments since the possibility for reduction of codimension is not assured in this setting.